\documentclass[12pt]{amsart}
\usepackage{fullpage,url,amssymb,enumerate,colonequals}

\usepackage{mathrsfs} 
\usepackage{MnSymbol}
\usepackage{extarrows}
\usepackage{lscape}
\usepackage[all,cmtip]{xy}

\usepackage[OT2,T1]{fontenc}

\usepackage{color}

\usepackage[
        colorlinks, citecolor=darkgreen,
        backref,
        pdfauthor={Nuno Freitas}, 
]{hyperref}

\usepackage{comment}

\newcommand{\F}{\mathbb{F}}

\newcommand{\Q}{\mathbb{Q}}

\newcommand{\Z}{\mathbb{Z}}
\newcommand{\Qbar}{{\overline{\Q}}}

\newcommand{\rhobar}{{\overline{\rho}}}

\newcommand{\fp}{\mathfrak{p}}

\DeclareMathOperator{\Gal}{Gal}

\newcommand{\vv}{\upsilon}

\newcommand{\SL}{\operatorname{SL}}

\numberwithin{equation}{section}

\newtheorem{theorem}{Theorem}

\theoremstyle{definition}

\theoremstyle{remark}

\definecolor{darkgreen}{rgb}{0,0.5,0}

\begin{document}

\title{An application of the symplectic argument to some Fermat-type Equations}

\author{Nuno Freitas}
\address{
University of British Columbia,
Department of Mathematics, 
Vancouver, BC V6T 1Z2
Canada }
\email{nunobfreitas@gmail.com}

\author{Alain Kraus}
\address{Universit\'e Pierre et Marie Curie - Paris 6,
Institut de Math\'ematiques de Jussieu,
4 Place Jussieu, 75005 Paris, 
France}
\email{alain.kraus@imj-prg.fr}


\begin{abstract}
  Let $p$ be a prime number. In the early 2000s, it was proved that the Fermat equations with coefficients 
  \[3x^p + 8y^p + 21z^p =0\quad \text{ and } \quad 3x^p + 4y^p + 5z^p=0 \] 
  do not admit non-trivial  solutions for a set of exponents $p$
  with Dirichlet density ${1/4}$ and ${1/8}$, respectively. 
  In this note, using a recent criterion to decide if 
  two elliptic curves over $\Q$ with certain types of 
  additive reduction at 2 have symplectically isomorphic
  $p$-torsion modules,  we improve these densities to  ${3/8}$.

  \bigskip
  
  \noindent {\sc{R\'esum\'e.}} {\bf Une application du crit\`ere symplectique \`a quelques \'equations de Fermat.} 
  Soit $p$ un nombre premier. Au d\'ebut des ann\'ees 2000, il a \'et\'e d\'emontr\'e que les \'equations de Fermat \`a coeficients 
  \[3x^p+8y^p+21z^p=0 \quad \hbox{et}\quad 3x^p+4y^p+5z^p=0\]
 ne poss\`edent pas de solutions non triviales pour un ensemble d'exposants $p$ de densit\'e de Dirichlet $1/4$ et $1/8$,  respectivement. Dans cette note, en utilisant un r\'esultat r\'ecent permettant  de d\'ecider si deux courbes elliptiques sur $\Q$, ayant un certain type de r\'eduction additive en $2$, ont leurs modules des points de $p$-torsion symplectiquement isomorphes, on am\'eliore ces densit\'es \`a $3/8$.

\end{abstract}

\maketitle

\section{Introduction}

In \cite{Serre87} J.P.-Serre raised questions regarding
specific instances of the \textit{Fermat equation with coefficients},
that is
\begin{equation}
  ax^p + by^p + cz^p=0 
  \label{E:FEWC}
\end{equation}
where $p$ is a prime and $a,b,c \in \Z$ are fixed and non-zero.

In \cite{HK2002} the second author and Halberstadt introduced the
{\it symplectic argument} as a complement to the modular method
and partly solved the questions raised by Serre,  together with other instances of the equation above.
Indeed, for some explicit choices of $a$, $b$ and $c$
they used Theorem~\ref{T:crit2} below and \cite[Lemma~1.7]{HK2002} to
show that \eqref{E:FEWC} does not have solutions
when the exponent $p$ belongs to certain congruence classes.
Another diophantine application of the symplectic argument and
Theorem~\ref{T:crit2} can be found in \cite{FS3}, where it was
used to solve the classical Fermat equation over $\Q(\sqrt{17})$
for a set of exponents with density 1/2.

We remark that the reason why the symplectic argument is necessary
is invisible in the proof of Fermat's Last Theorem where the
modular method has its origin. Indeed, after applying modularity and
level lowering results, one gets an isomorphism
\[
 \rhobar_{E,p} \sim \rhobar_{f,p}
\]
between the mod~$p$ representations attached to the Frey curve and
some newform $f$ with weight 2 and `small' level $N$. In the proof of FLT we
have $N=2$ and there are no candidate newforms $f$, giving a contradiction. 
In essentially every other application of the modular method there are candidates for $f$,
therefore more work is needed to obtain a contradiction to the previous isomorphism.
The symplectic argument is a tool that allows to obtain the desired contradiction
in certain cases. In particular, in the recent work \cite{F33p} the first author proved a 
new symplectic criterion (Theorem~\ref{T:crit1} below) and used it to solve the Generalized Fermat equation $x^3 + y^3 = z^p$ when $(-3/p) = -1$.

The purpose of this note is to further illustrate the strength of the symplectic argument, 
by combining the new and old criteria to improve two 
results originally obtained in \cite{HK2002}. More precisely, 
we will establish the following theorems.

\begin{theorem} Let $p > 7$ be a prime satisfying 
\[
p \equiv 5 \pmod{8} \qquad \text{ or } \qquad p \equiv 23 \pmod{24}.
\]
Then the Fermat equation 
\begin{equation}
3x^p + 8y^p + 21z^p =0          
\label{E:eq1}
\end{equation}
has no solutions $(x,y,z) \ne (0,0,0)$.
\label{T:main1}
\end{theorem}

\begin{theorem} Let $p \geq 5$ be a prime satisfying
\[
 p \equiv 5 \pmod{8} \qquad \text{ or } \qquad p \equiv 19 \pmod{24}.
\]
Then the Fermat equation 
\begin{equation}
3x^p + 4y^p+5z^p=0
\label{E:eq2}
\end{equation}
has no solutions $(x,y,z) \ne (0,0,0)$.
\label{T:main2}
\end{theorem}

 Let us mention that the coefficients $a$, $b$, $c$ of the Fermat equations \eqref{E:eq1} and \eqref{E:eq2} 
 do not satisfy non-trivial linear relations with coefficients in $\big\lbrace -1,0,1\big\rbrace$.
Conjecturally, we expect that for any prime number $p$ large enough, the two equations have at least one local obstruction i.e. there is at least a prime number $\ell$ such that they have no solutions over $\Q_{\ell}$ (see \cite[Conjecture  (C)]{HK2002}). Indeed, such is the case for $11\leq p<10^5$. For $p=3$, the curve of equation $3x^p+8x^p+21z^p=0$ has no points over  $\Q_3$ and  $\Q_7$; 
for $p=5$ and $p=7$, it has no local obstructions. The curve of equation  $3x^p+4y^p+5z^p=0$ has no local obstruction for $p=3$; for $p=5$ it has has no points over $\Q_{11}$ and for $p=7$  it  has no points over
$\Q_{29}$  and  $\Q_{43}$. 
Moreover, if $k$ is a fixed positive even integer, for any $p$ large enough  such that $q=kp+1$ is prime, the curves have a local obstruction at $q$ (see Proposition~3.3 in {\it loc. cit.}). For example, with $k=2$ such is the case as soon as $p\geq 11$.

\section{Two Symplectic criteria}

In this section we state the criteria we will use in the proofs. We first recall some
terminology.

Let $E$, $E'$ be elliptic curves over $\Q$ with $p$-torsion modules $E[p]$, $E'[p]$ and
Weil pairings $e_{E,p}$, $e_{E',p}$, respectively.
Let $\phi : E[p] \to E'[p]$ be an isomorphism of $G_\Q$-modules.
Then there is an element $r(\phi) \in \F_p^\times$ such that
\[ e_{E',p}(\phi(P), \phi(Q)) = e_{E,p}(P, Q)^{r(\phi)} \quad \text{for all $P, Q \in E[p]$.} \]
Note that for any $a \in \F_p^\times$ we have $r(a\phi) = a^2 r(\phi)$.
We say that $\phi$ is a \textit{symplectic isomorphism} if $r(\phi) = 1$ or,
more generally, $r(\phi)$ is a square in~$\F_p^\times$.
Fix a non-square~$r_p \in \F_p^\times$.
We say that $\phi$ is a \textit{anti-symplectic isomorphism} if $r(\phi) = r_p$
or, more generally, $r(\phi)$ is a non-square in~$\F_p^\times$.
Finally, we say that $E[p]$ and~$E'[p]$ are
\emph{symplectically} (or \emph{anti-symplectically}) \emph{isomorphic},
if there exists a symplectic (or anti-symplectic) isomorphism of~$G_\Q$-modules between them.

The following is \cite[Theorem~3]{F33p}. 

\begin{theorem}
  Let $E/\Q_2$ and $E'/\Q_2$ be elliptic curves with potentially good reduction. 
  Write $L =\Q_2^{\text{un}}(E[p])$ and $L' =\Q_2^{\text{un}}(E'[p])$.
  Write $\Delta_m(E)$ and $\Delta_m(E')$ for the minimal discriminant of $E$ and $E'$ respectively.
  Let $I_2 \subset \Gal(\Qbar_2 / \Q_2)$ be the inertia group.
  
  Suppose that $L = L'$ and $\Gal(L/\Q_2^{\text{un}}) \simeq \SL_2(\F_3)$.
  Then, $E[p]$ and $E'[p]$ are isomorphic $I_2$-modules for all prime $p \geq 3$. Moreover,
\begin{enumerate}
 \item if $(2/p)=1$ then $E[p]$ and $E'[p]$ are symplectically isomorphic $I_2$-modules. 
 \item if $(2/p)=-1$ then $E[p]$ and $E'[p]$ are symplectically isomorphic $I_2$-modules
 if and only if $\upsilon_2(\Delta_m(E)) \equiv \upsilon_2(\Delta_m(E')) \pmod{3}$. 
\end{enumerate}
Furthermore, $E[p]$ and $E'[p]$ cannot be both symplectic and anti-symplectic isomorphic $I_2$-modules.
\label{T:crit1}
\end{theorem}

The following is \cite[Proposition~2]{KO}.

\begin{theorem}
  Let $E$, $E'$ be elliptic curves over $\Q$ with minimal discriminants
  $\Delta$, $\Delta'$. Let $p$ be a prime such that $\rhobar_{E,p} \simeq \rhobar_{E',p}$.
  Suppose that $E$ and~$E'$ have multiplicative reduction at a prime $\ell \ne p$
  and that $p \nmid v_{\ell}(\Delta)$. Then $p \nmid v_{\ell}(\Delta')$,
  and the representations $\rhobar_{E,p}$ and $\rhobar_{E',p}$ are symplectically isomorpic
  if and only if $v_\ell(\Delta)/v_\ell(\Delta')$ is a square mod~$p$.
  \label{T:crit2}
\end{theorem}

\section{Proof of Theorem~\ref{T:main1}}

Suppose $(x,y,z)$ is a non-trivial primitive solution to \eqref{E:eq1}.
From \cite[Example~2.5]{HK2002} we know that the Frey curve $E_{x,y,z}$ attached to $(x,y,z)$
has minimal discriminant $\Delta_{x,y,z}$ given by
\begin{equation*}
\Delta_{x,y,z} = \begin{cases}
2^{10} \cdot 3^{2p-2} \cdot 7^2 \cdot (xyz)^{2p} & \text{ if } y \text{ is odd} \\
2^{-2} \cdot 3^{2p-2} \cdot 7^2 \cdot (xyz)^{2p} & \text{ if } y \text{ is even.} \
\end{cases}
\end{equation*}
Moreover, after applying the now classical modularity, irreducibilty and level lowering results
over $\Q$ we conclude that
\[
 \rhobar_{E_{x,y,z},p} \sim \rhobar_{f,\fp}
\]
where $f$ is a newform for $\Gamma_0(N)$ and weight 2 with level $N$ given by
\begin{equation*}
N = \begin{cases}
168 & \text{ if } y \text{ is odd} \\
42  & \text{ if } y \text{ is even.} \
\end{cases}
\end{equation*}
There is only one such newform at level 42 and two of them
at level 168. The three have rational coefficients hence correspond to isogeny classes
of elliptic curves. 
We note that the curves in the isogeny class with Cremona label `42a' have multiplicative reduction at $2$
while the curves of conductor 168 have potentially good reduction at 2. Furthermore,
their minimal extension $L/\Q_2^{un}$ of good reduction satisfies 
$\Gal(L/\Q_2^{un}) \simeq \SL_2(\F_3)$. 

We now divide the proof into two natural cases.

\noindent{Case I:} Suppose $y$ is even. Thus $E_{x,y,z}[p] \simeq E[p]$, where
$E=42a1$. It is proved in \cite[Example~2.5]{HK2002} that we get a contradiction 
with $(-2/p) = -1$.

\noindent{Case II:} Suppose $y$ is odd. There is an isomorphism 
$\phi : E_{x,y,z}[p] \simeq E[p]$, where
\[
 E = 168a1, \quad \Delta_E = 2^4 \cdot 3 \cdot 7 \qquad \text{ or } 
 \qquad E = 168b1, \quad \Delta_E = -2^4 \cdot 3^3 \cdot 7^4.
\]

Note that
\[
\vv_2(\Delta_{x,y,z}) = 10, \quad \vv_3(\Delta_{x,y,z}) \equiv -2, \quad \vv_7(\Delta_{x,y,z}) \equiv 2, 
\]
where the congruences are mod~$p$.

Suppose $(2/p) = -1$ and $E=\text{168a1}$. 
It follows from Theorem~\ref{T:main1} that $\phi$ is symplectic.
Thus Theorem~\ref{T:crit2} implies that 
$\vv_7(\Delta_{x,y,z})/\vv_7(\Delta_E) \equiv 2$ is a square mod~$p$, a contradiction.

Suppose $(2/p) = -1$ and $E=\text{168b1}$. 
It follows from Theorem~\ref{T:main1} that $\phi$ is symplectic.
Thus Theorem~\ref{T:crit2} implies that 
$\vv_7(\Delta_{x,y,z})/\vv_7(\Delta_E) \equiv 2/4$ is a square mod~$p$, a contradiction.

Suppose $(2/p) = 1$ and $E=\text{168a1}$. 
It follows from Theorem~\ref{T:main1} that $\phi$ is symplectic.
Thus Theorem~\ref{T:crit2} implies that 
$\vv_3(\Delta_{x,y,z})/\vv_3(\Delta_E) \equiv -2$ is a square mod~$p$.
This implies $(-1/p) = 1$.

Suppose $(2/p) = 1$ and $E=\text{168b1}$. 
It follows from Theorem~\ref{T:main1} that $\phi$ is symplectic.
Thus Theorem~\ref{T:crit2} implies that 
$\vv_3(\Delta_{x,y,z})/\vv_3(\Delta_E) \equiv -2/3$ is a square mod~$p$.
This implies $(-3/p) = 1$.

We therefore obtain a contradiction for all $y$ if one of the
following holds
\begin{itemize}
 \item $(-2/p) = -1$ and $(2/p) = -1$ or,
 \item $(-2/p) = -1$ and $(2/p) = (3/p) = 1$
\end{itemize}
which represents the  set of primes,  with density $3/8$, in the statement of the Theorem 1.

\section{Proof of Theorem~\ref{T:main2}}

Suppose $(x,y,z)$ is a non-trivial primitive solution to \eqref{E:eq2}.
From \cite[Proposition~2.3]{HK2002} and its proof 
we know that the Frey curve $E_{x,y,z}$ attached to $(x,y,z)$
has minimal discriminant $\Delta_{x,y,z}$ given by
\begin{equation*}
\Delta_{x,y,z} = \begin{cases}
2^{8} \cdot 3^{2} \cdot 5^2 \cdot (xyz)^{2p} & \text{ if } y \text{ is odd} \\
2^{-4} \cdot 3^{2} \cdot 5^2 \cdot (xyz)^{2p} & \text{ if } y \text{ is even.} \
\end{cases}
\end{equation*}
Moreover, after applying the now classical modularity, irreducibility and level lowering results
over $\Q$ we conclude that
$\rhobar_{E_{x,y,z},p} \sim \rhobar_{f,\fp}$
where $f$ is a newform for $\Gamma_0(N)$ and weight 2 with level $N$ given by
\begin{equation*}
N = \begin{cases}
120 & \text{ if } y \text{ is odd} \\
30  & \text{ if } y \text{ is even.} \
\end{cases}
\end{equation*}
There is only one such newform at level 30 and two of them
at level 120, each corresponding to an isogeny class of elliptic curves. 
Note that the curves in the isogeny class with Cremona label `30a' have multiplicative reduction at $2$
while the curves of conductor 120 have potentially good reduction at 2. Furthermore,
their minimal extension $L/\Q_2^{un}$ of good reduction satisfies 
$\Gal(L/\Q_2^{un}) \simeq \SL_2(\F_3)$. 

We now divide the proof into two natural cases.

\noindent {\sc Case I:} Suppose $y$ is even. Thus $E_{x,y,z}[p] \simeq E[p]$, where
\[
 E = 30a1, \qquad \Delta_E = -2^4 3^3 5^2.
\]
From Theorem~\ref{T:main2} the integers
\begin{itemize}
 \item $\vv_2(\Delta_{x,y,z})\vv_3(\Delta_{x,y,z})$ and $\vv_2(\Delta_E)\vv_3(\Delta_E)$,
 \item $\vv_2(\Delta_{x,y,z})\vv_5(\Delta_{x,y,z})$ and $\vv_2(\Delta_E)\vv_5(\Delta_E)$,
 \item $\vv_3(\Delta_{x,y,z})\vv_5(\Delta_{x,y,z})$ and $\vv_3(\Delta_E)\vv_5(\Delta_E)$,
\end{itemize}
must differ by multiplication by a square mod~$p$. This gives a contradiction with
\[
 (-2/p) = -1 \qquad \text{ or } \qquad (3/p) = -1.
\]
\noindent {\sc Case II:} Suppose $y$ is odd. Thus $\phi : E_{x,y,z}[p] \simeq E[p]$, where
\[
 E = 120a1, \quad \Delta_E = 2^4 \cdot 3^2 \cdot 5 \qquad \text{ or } 
 \qquad E = 120b1, \quad \Delta_E = -2^8 \cdot 3 \cdot 5.
\]

Note that $\vv_2(\Delta_{x,y,z}) = 8$ and 
$\vv_3(\Delta_{x,y,z}) \equiv \vv_5(\Delta_{x,y,z}) \equiv 2$ mod~$p$.

Suppose $(2/p) = -1$ and $E=\text{120a1}$. 
It follows from Theorem~\ref{T:main1} that $\phi$ is anti-symplectic.
Thus Theorem~\ref{T:crit2} implies that 
$\vv_3(\Delta_{x,y,z})/\vv_3(\Delta_E) \equiv 1$ is not a square mod~$p$, a contradiction.

Suppose $(2/p) = -1$ and $E=\text{120b1}$. 
It follows from Theorem~\ref{T:main1} that $\phi$ is symplectic.
Thus Theorem~\ref{T:crit2} implies that 
$\vv_5(\Delta_{x,y,z})/\vv_3(\Delta_E) \equiv 2$ is a square mod~$p$, a contradiction.

For the case $(2/p) = 1$ we cannot find further restrictions. 

We therefore obtain a contradiction for all $y$
if one of the following holds
\begin{itemize}
 \item $(2/p) = -1$ and $(-2/p) = -1$,
 \item $(2/p) = -1$ and $(3/p) = -1$.
\end{itemize}
The condition $(2/p) = -1$ means $p \equiv 3,8$ mod~$5$; if $p \equiv 5$ mod~8 
we have $(-2/p)= -1$ and the result follows in this case. Suppose $p \equiv 3$
mod~8, hence $(-2/p)= 1$. Now the condition $(3/p)=-1$ implies $p \equiv 1$ mod~3.
We get $p \equiv 19$~mod 24, as desired.  So we can conclude for a set of primes with density $3/8$.

\end{document}